\documentclass{amsart}
\usepackage{amsfonts,amssymb,amsmath,amsthm}
\usepackage{txfonts}
\usepackage{lmodern}
\usepackage{url}
\usepackage{enumerate}
\usepackage{graphicx}

\newtheorem{thrm}{Theorem}[section]
\newtheorem{lem}[thrm]{Lemma}
\newtheorem{prop}[thrm]{Proposition}
\newtheorem{cor}[thrm]{Corollary}
\theoremstyle{definition}
\newtheorem{definition}[thrm]{Definition}

\newtheorem{question}[thrm]{Question}
\newtheorem*{acknowledgement}{Acknowledgements}
\numberwithin{equation}{section}

\newcommand{\contract}{\mathord{\varparallelinv}}

\author{Lionel Pournin}
\address{LIPN, Universit{\'e} Paris 13, Villetaneuse, France}
\email{lionel.pournin@univ-paris13.fr}
\thanks{The author is partially supported by Ville de Paris {\'E}mergences project ``Combinatoire {\`a} Paris'' and by the ANR project SoS (Structures on Surfaces) number ANR-17-CE40-0033.}

\begin{document}

\title[Eccentricities in flip-graphs]{Eccentricities in the flip-graphs of~convex~polygons}

\begin{abstract}
The flip-graph of a convex polygon $\pi$ is the graph whose vertices are the triangulations of $\pi$ and whose edges correspond to flips between them. The eccentricity of a triangulation $T$ of $\pi$ is the largest possible distance in this graph from $T$ to any triangulation of $\pi$. It is well known that, when all $n-3$ interior edges of $T$ are incident to the same vertex, the eccentricity of $T$ in the flip-graph of $\pi$ is exactly $n-3$, where $n$ denotes the number of vertices of $\pi$. Here, this statement is generalized to arbitrary triangulations. Denoting by $n-3-k$ the largest number of interior edges of $T$ incident to a vertex, it is shown that the eccentricity of $T$ in the flip-graph of $\pi$ is exactly $n-3+k$, provided $k\leq{n/2-2}$. Inversely, the eccentricity of a triangulation, when small enough, allows to recover the value of $k$. More precisely, if $k\leq{n/8-5/2}$, it is also shown that $T$ has eccentricity $n-3+k$ if and only if exactly $n-3-k$ of its interior edges are incident to a given vertex. When $k>n/2-2$, bounds on the eccentricity of $T$ are also given and discussed.
\end{abstract}
\maketitle

\section{Introduction}
A triangulation of a convex polygon is a set of non-crossing line segments, the edges of the triangulation, whose endpoints are vertices of the polygon. In order for these segments to decompose the polygon into triangles, a triangulation is also required to be maximal for the inclusion. One can put a nice graph structure on the set of triangulations of a convex polygon by considering local operations called \emph{flips}. A flip consists in removing an edge from a triangulation, provided it is not also an edge of the polygon, and then replacing it by the only other line segment such that the resulting object is still a triangulation of the same polygon. The \emph{flip-graph} of a polygon is then the graph whose vertices are the triangulations of that polygon and whose edges correspond to flips. Flip-graphs turn up in a wide variety of topics including, for instance, discrete geometry \cite{DeLoeraRambauSantos2010,Santos2000,Santos2005}, geometric topology \cite{BellDisarloTang2018,DisarloParlier2018,ParlierPournin2017}, probability \cite{Budzinski2017,CaputoMartinelliSinclairStauffer2015}, computer science \cite{BoseHurtado2009,CulikWood1982}, or biology \cite{SempleSteel2003}. 

The flip-graph of a convex polygon has a number of remarkable properties, including that of being the graph of a polytope, the associahedron \cite{Lee1989}. It has been a long-standing open problem to find its diameter for polygons with any number of vertices \cite{SleatorTarjanThurston1988}. While this particular problem is now solved \cite{Pournin2014}, the geometry of this and related flip-graphs is still not fully understood. For instance, computing the distances in these graphs is instrumental for a number of applications \cite{SempleSteel2003,SleatorTarjanThurston1988}. While such computations are proven to be hard in more general settings \cite{AichholzerMulzerPilz2015,LubiwPathak2015,Pilz2014}, there is as yet no such result, or a polynomial algorithm to compute distances in the flip-graph of a convex polygon. 
Halfway between the notions of diameter and distance mentioned in these problems, one finds eccentricity.

\begin{definition}
Consider a graph $G$ and a vertex $v$ of $G$. The \emph{eccentricity} of $v$ in $G$ is the largest possible distance in $G$ from $v$ to any vertex of $G$.
\end{definition}

In order to estimate the diameter of the flip-graph of a convex polygon $\pi$, the authors of \cite{SleatorTarjanThurston1988} explicitly build a path in this graph from an arbitrary triangulation $T$ of $\pi$ to a \emph{comb}, a triangulation whose all interior edges are incident to a given vertex of $\pi$. In particular, they prove that $T$ can be transformed into a comb by a sequence of at most $n-3$ flips, where $n$ denotes the number of vertices of $\pi$. It turns out that this transformation can require exactly $n-3$ flips. In other words, the eccentricity of a comb in the flip-graph of $\pi$ is $n-3$. It is interesting to note that $n-3$ is precisely the number of interior edges of any triangulation of a convex polygon with $n$ vertices. The above statement on the eccentricity of combs is generalized here to arbitrary triangulations.

The first main result in this paper is the following theorem.

\begin{thrm}\label{Thm.Main.1}
Consider a convex polygon $\pi$ with $n$ vertices. If a vertex of $\pi$ is incident to exactly $n-3-k$ interior edges of a triangulation $T$ of $\pi$, where $k\leq{n/2-2}$, then $T$ has eccentricity $n-3+k$ in the flip-graph of $\pi$.
\end{thrm}

Inversely, the eccentricity of a triangulation in the flip-graph of $\pi$, when small enough, allows to recover the number of its edges incident to one of its vertices. More precisely, the following theorem is obtained here as well.

\begin{thrm}\label{Thm.Main.2}
Consider a convex polygon $\pi$ with $n$ vertices. If $k\leq{n/8-5/2}$, then a triangulation of $\pi$ has eccentricity $n-3+k$ in the flip-graph of $\pi$ if and only if exactly $n-3-k$ of its interior edges are incident to a vertex.
\end{thrm}

In order to prove these theorems, two lower bounds on the eccentricity of the triangulations of a convex polygon will be established. The first bound, proven in Section \ref{Sec.close} along with Theorem \ref{Thm.Main.1}, is sharp and it is only valid for triangulations with a vertex incident to more than half of their interior edges. Note that these triangulations are close to combs. The second bound, obtained in Section \ref{Sec.away} and used to prove Theorem \ref{Thm.Main.2}, is not sharp but it holds in general. In particular, it gives new information on the eccentricities of triangulations far away from combs. In Section \ref{Sec.prelim}, the notions and tools from \cite{Pournin2014,SleatorTarjanThurston1988} that will be used to obtain these lower bounds on the eccentricities are described. Additional results and remarks on the behavior of the eccentricities far away from combs are given in Section \ref{Sec.conc}. In particular, it will be shown that Theorems \ref{Thm.Main.1} and~\ref{Thm.Main.2} do not extend to arbitrary values of $k$. Section \ref{Sec.conc} ends with a couple of questions.

\section{Preliminary notions and tools}\label{Sec.prelim}

Let $\pi$ be a convex polygon. In the following, an \emph{edge on $\pi$} is a set containing exactly two vertices of $\pi$. Using this terminology, the edges \emph{of} $\pi$ are the edges on $\pi$ whose convex hull is disjoint from the interior of $\pi$. Two edges on $\pi$ are crossing when their convex hulls have non-disjoint interiors. A \emph{triangulation of $\pi$} is a set of pairwise non-crossing edges on $\pi$ that is maximal for the inclusion. Note that all the edges of $\pi$ are contained in any of its triangulations, and they will be referred to as the \emph{boundary edges} of these triangulations. The other edges of a triangulation will be called its \emph{interior edges}. The \emph{flip-graph} of $\pi$ is a graph whose vertices are the triangulations of $\pi$. Two triangulations are adjacent in this graph when they differ by a single edge, or equivalently, when they are related by a flip. 

In the following, a \emph{path of length $k$} between two triangulations $T$ and $U$ of a convex polygon $\pi$ is a sequence of $k$ flips that transform one of these triangulations into the other. A shortest path between $T$ and $U$ will also be called a \emph{geodesic}, and the length of any geodesic between $T$ and $U$ will be denoted by $d(T,U)$. The following straightforward proposition is used, sometimes implicitly, in a number of papers about the geometry of flip-graphs (see \cite{SleatorTarjanThurston1988}).

\begin{prop}\label{Intro.prop.1}
Let $T$ be a triangulation of a convex polygon $\pi$ with $n$ vertices. If some vertex of $\pi$ is incident to exactly $n-3-k$ interior edges of $T$ then, for any triangulation $U$ of $\pi$, $d(T,U)\leq{n-3+k}$.
\end{prop}

As in \cite{ParlierPournin2017,Pournin2014,SleatorTarjanThurston1988}, proving lower bounds on the distances (and in the case at hand, on the eccentricities) in flip-graphs is the hard part. In order to do so, some of the techniques developed in \cite{Pournin2014} are borrowed here.

If $a$ and $b$ are two vertices of a convex polygon $\pi$ such that $b$ immediately follows $a$ clockwise, the pair $(a,b)$ will be called a \emph{clockwise-oriented boundary edge of $\pi$}. If $\pi$ has at least four vertices, then removing the edge $\{a,b\}$ from a triangulation $T$ of $\pi$, and replacing $a$ by $b$ in all the remaining edges results in a triangulation of a smaller polygon (see Proposition 2 in \cite{Pournin2014}). This operation will be referred to as the \emph{deletion of vertex $a$ from $T$}, and the resulting triangulation will be denoted by $T\contract{a}$. Note in particular that $T\contract{a}$ is a triangulation of the polygon whose vertex set is obtained by removing $a$ from the vertex set of $\pi$.

A flip performed within $T$ will be called \emph{incident to $\{a,b\}$} if it affects the triangle of $T$ incident to $\{a,b\}$. The following two lemmas are proven in \cite{Pournin2014}.

\begin{lem}\label{Lem.intro.0}
Consider two triangulations $T$ and $U$ of a convex polygon $\pi$ with at least four vertices and a clockwise-oriented boundary edge $(a,b)$ of $\pi$. If $f$ flips are incident to $\{a,b\}$ along a geodesic between $T$ and $U$, then
$$
d(T,U)\geq{d(T\contract{a},U\contract{a})+f}\mbox{.}
$$
\end{lem}

Observe that a triangulation $T$ of $\pi$ decomposes $\pi$ into a set of triangles. If such a triangle shares two of its edges $\{a,b\}$ and $\{a,c\}$ with $\pi$, it will be called \emph{an ear of $T$ in $a$}. Equivalently, $T$ has an ear in $a$ when none of its interior edges is incident to $a$. Using this notion, a second lemma can be stated.

\begin{lem}\label{Lem.intro.1}
Consider a clockwise-oriented boundary edge $(a,b)$ of a convex polygon $\pi$ with at least four vertices. If $T$ is a triangulation of $\pi$ with an ear in $b$ and if $U$ is a triangulation of $\pi$ with at least two interior edges incident to $b$, then there exists an $x\in\{a,b\}$ such that $d(T,U)\geq{d(T\contract{x},U\contract{x})+2}$.
\end{lem}

\section{Eccentricities close to combs}\label{Sec.close}

In order to bound the eccentricity of a triangulation close to combs, the following notion of a shelling at a vertex will be used.

\begin{definition}\label{Def.shelling}
Consider a triangulation $T$ of a polygon $\pi$ with $n$ vertices and a vertex $v$ of $\pi$ incident to exactly $n-3-k$ interior edges of $T$. Consider the vertices of $\pi$ that are not adjacent to $v$ by an edge of $T$. A \emph{shelling of $T$ at $v$} is an ordering $a_1$, ..., $a_k$ of these vertices such that for all $i\in\{1, ..., k\}$, the edges of $T$ not incident to any of the vertices $a_i$, ..., $a_k$ still form a triangulation of a convex polygon. 
\end{definition}

\begin{figure}
\begin{centering}
\includegraphics{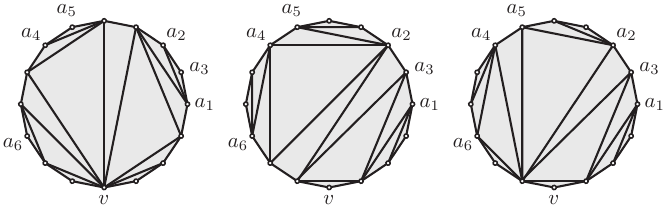}
\caption{A shelling of a triangulation $T$ at a vertex $v$ (left) and two triangulations compatible with that shelling.}\label{Close.fig.1}
\end{centering}
\end{figure}

An example of a shelling at a vertex is depicted on the left of Fig.~\ref{Close.fig.1}. Consider the shelling in this figure, and observe that $T$ has an ear in $a_6$. Further note that $a_1$, ..., $a_5$ is a shelling at $v$ of the triangulation obtained by removing that ear from $T$. This is a general property: a triangulation $T$ always has an ear in the largest-indexed vertex $a_k$ of any of its shellings at a given vertex, and removing $a_k$ from the shelling results in a shelling of the triangulation obtained by removing the ear in $a_k$ from $T$. This property allows for inductive proofs.

The purpose of the next definition is to provide a family of triangulations of a convex polygon $\pi$ that, as will be shown later in the section, achieve the largest distance in the flip-graph of $\pi$ to a given triangulation of $\pi$.

\begin{definition}\label{Def.close.1}
Let $T$ be a triangulation of a polygon $\pi$ with $n$ vertices and $v$ a vertex of $\pi$ incident to exactly $n-3-k$ interior edges of $T$. A triangulation $U$ of $\pi$ will be called \emph{compatible with a shelling} $a_1$, ..., $a_k$ of $T$ at $v$ when $U$ has an ear in $v$ and, for any $i\in\{1, ..., k\}$, $a_i$ is incident to at least two interior edges of $U$ whose other vertex is not among $a_i$, ..., $a_k$.
\end{definition}

Consider the triangulation $T$ and its shelling at $v$ shown on the left of Fig.~\ref{Close.fig.1}. The two other triangulations depicted in the same figure are examples of triangulations compatible with that shelling. Note that two interior edges of the triangulation in the center are incident to $a_2$ and to a vertex $a_i$ such that $i>2$. This triangulation is compatible with the considered shelling nonetheless.
Indeed, $a_2$ is also incident to two interior edges whose other vertex is adjacent to $v$ by an edge of $T$.

According to the following lemma, the notion of compatibility with a shelling at a vertex is preserved under well-chosen vertex deletions.

\begin{lem}\label{Lem.sec1.1}
Let $T$ be a triangulation of a polygon $\pi$ with at least four vertices, and $v$ a vertex of $\pi$ incident to exactly $n-3-k$ interior edges of $T$. Consider a triangulation $U$ compatible with a shelling $a_1$, ..., $a_k$ of $T$ at $v$. If $x$ is equal to $a_k$ or to the vertex of $\pi$ that immediately precedes $a_k$ clockwise, then $U\contract{x}$ is compatible with at least one shelling of $T\contract{x}$ at $v$.
\end{lem}

\begin{proof}
Let $x$ be equal to $a_k$ or to the vertex of $\pi$ that immediately precedes $a_k$ clockwise. By the definition of a shelling, $T$ has an ear in $a_k$, and this ear is removed when $x$ is deleted from $T$. In particular, if $x$ is distinct from vertices $a_1$ to $a_{k-1}$, then $a_1$, ..., $a_{k-1}$ is a shelling of $T\contract{x}$ at $v$. If $x$ is equal to $a_j$, where $1\leq{j}<k$, then it will be assumed that $a_k$ is relabeled $a_j$ as a vertex of $T\contract{x}$ and $U\contract{x}$. In this case, $a_1$, ..., $a_{k-1}$ is still a shelling of $T\contract{x}$ at $v$. In the remainder of the proof, it is shown that $U\contract{x}$ is compatible with that shelling.

Denote by $y$ the vertex of $\pi$ that follows $x$ clockwise. Consider the triangle of $U$ incident to $\{x,y\}$ and call $z$ its third vertex. Observe that $z$ is the only vertex of $\pi$ distinct from $x$ such that the number of interior edges incident to $z$ is less in $U\contract{x}$ than in $U$. Hence, if $z$ is distinct from vertices $a_1$ to $a_{k-1}$, then $U\contract{x}$ is necessarily compatible with the shelling $a_1$, ..., $a_{k-1}$ of $T\contract{x}$ at $v$. If $z$ is equal to $a_i$, where $1\leq{i}<k$, then by Definition \ref{Def.close.1}, it is incident to at least two interior edges of $U$ whose other vertex is not among $a_i$, ..., $a_k$. If the other vertex of such an edge is not equal to $x$, then this edge is still an interior edge of $U\contract{x}$ incident to $z$ whose other vertex is not among $a_i$, ..., $a_{k-1}$. It is possible, though, that one of the edges incident to $z$ whose other vertex is not among $a_i$, ..., $a_k$ is $\{x,z\}$. If this happens, a new such interior edge needs to be found in $U\contract{x}$. If $a_k$ has not been relabeled $a_j$, then this edge will be $\{a_k,z\}$ (note in particular, that $a_k$ is not among $a_i$, ..., $a_{k-1}$). If $a_k$ has been relabeled $a_j$, then this edge will be $\{a_j,z\}$ (in this case $x=a_j$ before the deletion and it follows that $a_j$ is not among $a_i$, ..., $a_{k-1}$ before or after the deletion). 
As a consequence, $z$ is still incident to at least two interior edges of $U\contract{x}$ whose other vertex is not among $a_i$, ..., $a_{k-1}$. This proves that $U\contract{x}$ is compatible with the shelling $a_1$, ..., $a_{k-1}$ of $T\contract{x}$ at $v$.
\end{proof}

Consider a triangulation $T$ of a polygon with $n$ vertices and a vertex $v$ of $\pi$ incident to exactly $n-3-k$ interior edges of $T$. Note that $k>n/2-2$ if and only if at most half of the interior edges of $T$ are incident to $v$. It turns out that, in this case, no triangulation of $\pi$ can be compatible with a shelling of $T$ at $v$. Indeed, consider such a shelling $a_1$, ..., $a_k$ and assume that a triangulation $U$ is compatible with it. By definition, each $a_i$ is incident to at least two interior edges of $U$ whose other vertex is not among $a_i$, ..., $a_k$. In particular, two such incidences cannot be to the same edge. Hence, $2k$ is a lower bound on the number of interior edges of $U$ incident to $a_1$, ..., $a_k$. In addition, the interior edge of $U$ bounding the ear in $v$ is not incident to any of the vertices $a_1$ to $a_k$. Therefore, $2k$ is not greater than $n-4$, or equivalently $k\leq{n/2-2}$. Note that, as a consequence, the statement of the following theorem is void for these values of $k$.

\begin{thrm}\label{MPThm.1}
Let $T$ be a triangulation of a polygon $\pi$. If $v$ is a vertex of $\pi$ incident to exactly $n-3-k$ interior edges of $T$, then the distance between $T$ and any triangulation compatible with a shelling of $T$ at $v$ is at least $n-3+k$.
\end{thrm}
\begin{proof}
The theorem shall be proven by induction on $k$. If $k=0$, then all the interior edges of $T$ are incident to $v$. In this case, the only possible shelling of $T$ at $v$ is empty and any triangulation of $\pi$ with an ear in $v$ is compatible with this shelling. Since the triangulations of $\pi$ all have $n-3$ interior edges and since a flip removes a single edge, the desired statement holds.

Now assume that $k\geq1$. Consider a triangulation $U$ of $\pi$ compatible with a shelling $a_1$, ..., $a_k$ of $T$ at $v$. By the definition of a shelling, $T$ has an ear in $a_k$. Moreover, according to Definition \ref{Def.close.1}, $a_k$ is incident to at least two interior edges of $U$. Hence, by Lemma \ref{Lem.intro.1},
\begin{equation}\label{MPThm.1.eq.1}
d(T,U)\geq{d(T\contract{x},U\contract{x})+2}\mbox{,}
\end{equation}
where $x$ is equal to $a_k$ or to the vertex of $\pi$ that precedes $a_k$ clockwise. According to Lemma \ref{Lem.sec1.1}, $U\contract{x}$ is compatible with some shelling of $T\contract{x}$ at $v$. Moreover, $v$ is still incident to $n-3-k$ interior edges of $T\contract{x}$. Now observe that
$$
n-3-k=(n-1)-3-(k-1)\mbox{.}
$$

Hence, by induction,
$$
d(T\contract{x},U\contract{x})\geq{(n-1)-3+k-1}\mbox{.}
$$

Combining this inequality with (\ref{MPThm.1.eq.1}) completes the proof.
\end{proof}

As argued above, no triangulation can be compatible with a shelling of a triangulation $T$ at a vertex $v$ incident to at most half of the interior edges of $T$. If, on the contrary, $v$ is incident to more than half of the interior edges of $T$, then according to Lemma \ref{Lem.sec1.2}, such triangulations exist. In order to prove that lemma, the following proposition is needed, that can be thought of as a very simple combinatorial variant of the ham-sandwich theorem. In the proof of Lemma \ref{Lem.sec1.2}, this proposition will be invoked for the vertex set of a polygon with the natural clockwise ordering.

\begin{prop}\label{Prop.close.1}
Let $V$ be a totally ordered finite set. If a subset $S$ of $V$ has cardinality at most $(|V|-3)/2$, then there exists an $x\in{V\mathord{\setminus}S}$ such that there are exactly $i$ elements of $S$ and exactly $2i+1$ elements of $V$ less than $x$.
\end{prop}
\begin{proof}
First observe that, since $|S|\leq(|V|-3)/2$, then $V$ must contain at least three elements. For any $y\in{V}$, respectively call $f_V(y)$ and $f_S(y)$ the number of elements of $V$ less than $y$ and the number of elements of $S$ less than $y$. When $y$ is the smallest element of $V$, the sum $f_V(y)-2f_S(y)-1$ is negative. As $|S|\leq(|V|-3)/2$, this sum is positive when $y$ is the largest element of $V$. Further observe that this sum varies by at most one from an element of $V$ to the next. Hence there must exist a $x\in{V}$ such that $f_V(x)-2f_S(x)-1=0$. In other words, exactly $i$ elements of $S$ and exactly $2i+1$ elements of $V$ are less than $x$, where $i=f_S(x)$. Assuming that $x$ is the largest such element of $V$, then $f_V(y)-2f_S(y)-1$ must be positive for the element $y$ of $V$ following $x$, proving that $x$ cannot belong to $S$.
\end{proof}


A constructive proof can now be given that, for a triangulation $T$ and a vertex $v$ of $T$ incident to more than half of the interior edges of $T$, there always exists a triangulation compatible with any shelling of $T$ at $v$.

\begin{lem}\label{Lem.sec1.2}
Let $T$ be a triangulation of a convex polygon $\pi$ with $n$ vertices and $v$ a vertex of $\pi$ incident to $n-3-k$ interior edges of $T$. If $k\leq{n/2-2}$, then there exists a triangulation of $\pi$ compatible with any shelling of $T$ at $v$.
\end{lem}
\begin{proof}
Assume that $k\leq{n/2-2}$ and consider a shelling $a_1$, ..., $a_k$ of $T$ at $v$. Let $m$ be an integer such that $3\leq{m}\leq{n-1}$ and $\pi'$ a polygon whose vertex set is made up of $m$ vertices of $\pi$ distinct from $v$. Note that there is at least one such polygon. Indeed, $n$ must be at least $4$ because $0\leq{k}\leq{n/2-2}$.
Call $S(\pi')$ the set of the vertices among $a_1$ to $a_k$ such that the two edges of $\pi$ incident to these vertices are both also edges of $\pi'$. Assume that $S(\pi')$ has cardinality at most $(m-3)/2$.

It will be proven by induction on $m$ that there exists a triangulation $U'$ of $\pi'$ such that, if $a_i$ belongs to $S(\pi')$ then it is incident to at least two interior edges of $U'$ whose other vertex is not in $S(\pi')\cap\{a_i, ..., a_k\}$. Note that, when $m=n-1$, the elements of $S(\pi')$ are exactly $a_1$ to $a_k$. Hence, in this case, gluing the ear in $v$ to $U'$ results in a triangulation compatible with the considered shelling of $T$.

If $S(\pi')$ is empty then any triangulation of $\pi'$ has the desired property. Otherwise, say that $j$ is the smallest index such that $a_j$ belongs to $S(\pi')$. Order the set $V$ of the vertices of $\pi'$ distinct from $a_j$ clockwise from the one that follows $a_j$ to the one that precedes it. Note that
\begin{equation}\label{Close.Exist.eq.1}
|V\cap{S(\pi')}|\leq\frac{|V|-4}{2}\mbox{.}
\end{equation}

Proposition \ref{Prop.close.1} can therefore be invoked for the set $V$ ordered clockwise and its subset $S(\pi')\mathord{\setminus}\{a_j\}$. It provides a vertex $x$ in $V\mathord{\setminus}S(\pi')$ such that the set $V_1$ of the vertices in $V$ less than $x$ satisfies
\begin{equation}\label{Close.Exist.eq.2}
|V_1\cap{S(\pi')}|=\frac{|V_1|-1}{2}\mbox{.}
\end{equation}

Combining (\ref{Close.Exist.eq.1}) and (\ref{Close.Exist.eq.2}) yields:
\begin{equation}\label{Close.Exist.eq.3}
|[V\mathord{\setminus}V_1]\cap{S(\pi')}|\leq\frac{|V\mathord{\setminus}V_1|-3}{2}\mbox{.}
\end{equation}

Hence, invoking Proposition \ref{Prop.close.1} again, but with $V\mathord{\setminus}V_1$ provides a vertex $y$ in $[V\mathord{\setminus}V_1]\mathord{\setminus}S(\pi')$ such that the set $V_2$ of the vertices less than $y$ in $V\mathord{\setminus}V_1$ satisfies
\begin{equation}\label{Close.Exist.eq.4}
|V_2\cap{S(\pi')}|=\frac{|V_2|-1}{2}\mbox{.}
\end{equation}

Further call $V_3=V\mathord{\setminus}[V_1\cup{V_2}]$, and observe that, by (\ref{Close.Exist.eq.3}) and (\ref{Close.Exist.eq.4}),
\begin{equation}\label{Close.Exist.eq.5}
|V_3\cap{S(\pi')}|\leq\frac{|V_3|-2}{2}\mbox{.}
\end{equation}

Now denote by $\pi'_1$, $\pi'_2$, and $\pi'_2$ the polygons whose vertex sets are, respectively, $V_1\cup\{x,a_j\}$, $V_2\cup\{y,a_j\}$, and $V_3\cup\{a_j\}$.
\begin{figure}[b]
\begin{centering}
\includegraphics{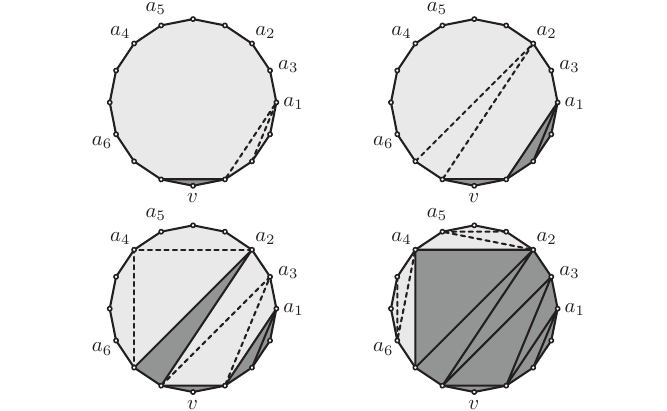}
\caption{An example of the construction in the proof of Lemma \ref{Lem.sec1.2}.
}\label{Close.fig.2}
\end{centering}
\end{figure}
These three polygons can be alternatively obtained by cutting $\pi'$ along the edges $\{a_j,x\}$ and $\{a_j,y\}$. The construction is illustrated in Fig. \ref{Close.fig.2}, where $\{a_j,x\}$ and $\{a_j,y\}$ are shown as dotted lines for all the polygons $\pi'$ (colored light grey) inductively considered in the case of the shelling from Fig. \ref{Close.fig.1}. Note that the triangulation resulting from the procedure, depicted on the bottom-right of Fig. \ref{Close.fig.2}, is the one in the center of Fig. \ref{Close.fig.1}.

Let $l\in\{1,2,3\}$. Observe that
$$
S(\pi'_l)=V_l\cap{S(\pi')}\mbox{.}
$$

Therefore, calling $m_l$ the number of vertices of $\pi_l$, it follows from either (\ref{Close.Exist.eq.2}), (\ref{Close.Exist.eq.4}), or (\ref{Close.Exist.eq.5}) that $|S(\pi'_l)|\leq(m_l-3)/2$. Hence, by induction, there exists a triangulation $U_l$ of $\pi'_l$ such that if $a_i$ belongs to $S(\pi'_l)$, then it is incident to at least two interior edges of $U'_l$ whose other vertex does not belong to $S(\pi'_l)\cap\{a_i, ..., a_k\}$. Now, call $U'=U'_1\cup{U'_2}\cup{U'_3}$. Consider a vertex among $a_1$ to $a_k$ that belongs to $S(\pi')$, say $a_i$. Recall that $S(\pi'_1)$, $S(\pi'_2)$, $S(\pi'_3)$, and $\{a_j\}$ form a partition of $S(\pi')$. Hence, if $i\neq{j}$, then there exists a unique $l\in\{1,2,3\}$ such that $a_i\in{S(\pi'_l)}$. In this case, as proven above, $a_i$ is incident to at least two interior edges of $U'_l$ whose other vertex is not in $S(\pi'_l)\cap\{a_i, ..., a_k\}$. Since $a_j$ is the only vertex of $\pi'_l$ in $S(\pi')\mathord{\setminus}S(\pi'_l)$ and since $j<i$, then $a_i$ is incident to at least two interior edges of $U'$ whose other vertex is not in $S(\pi')\cap\{a_i, ..., a_k\}$. If $i=j$ then, by construction, $a_i$ is incident to $\{a_i,x\}$ and $\{a_i,y\}$. As $x$ and $y$ do not belong to $S(\pi')$, the lemma is proven. 
\end{proof}

Theorem \ref{Thm.Main.1} can now be established.

\begin{proof}[Proof of Theorem \ref{Thm.Main.1}]
Consider a triangulation $T$ of a convex polygon $\pi$ with $n$ vertices. Further consider a vertex $v$ of $\pi$ incident to exactly $n-3-k$ interior edges of $T$. If $k\leq{n/2-2}$ then Lemma \ref{Lem.sec1.2} provides a triangulation of $\pi$ compatible with a shelling of $T$ at $v$. According to Theorem \ref{MPThm.1}, the eccentricity of $T$ in the flip-graph of $\pi$ is then at least $n-3+k$. By Proposition \ref{Intro.prop.1}, this is sharp.
\end{proof}

\section{Eccentricities away from combs}\label{Sec.away}

The main purpose of the section is to prove Theorem \ref{Thm.Main.2}. This will be done by establishing a lower bound on the eccentricities in the flip-graph of a convex polygon that holds for any triangulation regardless of its distance to a comb.

Consider a polygon $\pi$ with $n$ vertices. If $a$ and $b$ are any two vertices of $\pi$, the ordered pair $(a,b)$ will be called an \emph{oriented edge on $\pi$}. Such an edge splits the vertices of $\pi$ into two parts. More precisely, the vertices of $\pi$ can be labeled as $v_0$ to $v_{n-1}$ in such a way that $v_0=a$ and, for all $i\in\{1, ..., n-1\}$, $(v_{i-1},v_i)$ is a clockwise-oriented boundary edge of $\pi$. For some index $l$, that will be called the \emph{length of $(a,b)$}, vertices $b$ and $v_l$ coincide. If $l$ is not less than $2$, then the vertices $v_1$ to $v_{l-1}$ will be referred to as the vertices of $\pi$ \emph{on the left of $(a,b)$}. Similarly, if $l$ is not greater than $n-2$, then $v_{l+1}$ to $v_{n-1}$ will be the vertices of $\pi$ \emph{on the right of $(a,b)$}. Note that, while the lengths of $(a,b)$ and $(b,a)$ do not necessarily coincide, these two lengths always sum to $n$. As a result, for any vertex $c$ of $\pi$ on the right of $(a,b)$, the lengths of $(a,b)$, $(b,c)$, and $(c,a)$ also sum to $n$.

In the sequel, the following straightforward result will be needed.

\begin{prop}\label{Prop.far.1}
Let $T$ be a triangulation of a polygon $\pi$ and $(a,b)$ an oriented edge on $\pi$ such that $\{a,b\}$ is an edge of $T$. If $(a,b)$ has length at least $2$, then $T$ has an ear in a vertex of $\pi$ on the left of $(a,b)$.
\end{prop}
\begin{proof}
Assume that $(a,b)$ has length at least $2$. In this case, as $\{a,b\}$ is an edge of $T$, cutting $T$ along this edge results in a triangulation $U$ of the polygon placed on the left of $(a,b)$. 
If $(a,b)$ has length exactly $2$, then this polygon has three vertices and $U$ is made up of a single triangle. This triangle is an ear of $T$ in the only vertex of $\pi$ left of $(a,b)$. If $(a,b)$ has length at least $3$, then $U$ has at least two ears (see for instance \cite{Meisters1980}). At least one of these ears must be an ear in some vertex $v$ of $\pi$ on the left of $(a,b)$ and, therefore, it is also an ear of $T$ in $v$.
\end{proof}

Using these notions, one can define, for every triangulation $T$ of $\pi$, a set of triangulations of $\pi$ that will be at a reasonably large distance from $T$ in the flip-graph of $\pi$, provided $T$ is far away from every comb.

Consider an oriented edge $(a,b)$ on $\pi$, of length at most $\lceil{n/2}\rceil-1$ such that $\{a,b\}$ belongs to $T$. Call $\Omega(T,a,b)$ the set of all the triangulations $U$ of $\pi$ such that all the interior edges shared by $T$ and $U$ are incident to $a$ or to $b$, and all the vertices of $\pi$ on the left of $(a,b)$ are incident to at least two interior edges of $U$.

The structure of the triangulations in $\Omega(T,a,b)$ makes it possible to obtain reasonable lower bounds on their distance to $T$ by using Lemma \ref{Lem.intro.1} inductively.

\begin{thrm}\label{Thm.far.1} 
Consider a triangulation $T$ of a convex polygon $\pi$ with $n$ vertices. Let $(a,b)$ be an oriented edge on $\pi$ such that $\{a,b\}$ belongs to $T$. Denote the length of $(a,b)$ by $l$. Further consider a triangulation $U$ in $\Omega(T,a,b)$. Call $m$ the number of vertices of $\pi$ on the right of $(a,b)$ that are adjacent to $a$ or to $b$ by an interior edge of $T$ and to a vertex of $\pi$ that is not on the right of $(a,b)$ by an interior edge of $U$. The following inequality holds:
$$
d(T,U)\geq{n-m+l-5}\mbox{.}
$$
\end{thrm}

\begin{proof}
The theorem will be proven by induction on $l$. By the definition of $\Omega(T,a,b)$, the only possible edges common to $T$ and $U$ are incident to $a$ or to $b$. Moreover, at most $m$ vertices of $\pi$ are adjacent to $a$ or to $b$ by an interior edge common to $T$ and $U$, and these are necessarily on the right of $(a,b)$. As only one such vertex can be adjacent to both $a$ and $b$ in a triangulation, $T$ and $U$ share at most $m+1$ interior edges. Hence all the interior edges of $T$ must be flipped in order to obtain $U$, except possibly $m+1$ of them. As a consequence, $d(T,U)\geq{n-m-4}$. In particular, the desired result holds when $l$ is equal to $1$.

Now assume that $l\geq2$. In this case, by Proposition \ref{Prop.far.1}, $T$ has an ear in some vertex $c$ of $\pi$ on the left of $(a,b)$. Moreover, according to the definition of $\Omega(T,a,b)$, at least two interior edges of $U$ are incident to $c$. Hence, Lemma \ref{Lem.intro.1} yields:
\begin{equation}\label{Thm.far.1.eq.1}
d(T,U)\geq{d(T\contract{x},U\contract{x})+2}\mbox{,}
\end{equation}
where $x$ is either equal to $c$ or to the vertex of $\pi$ that immediately precedes $c$ clockwise. If $x$ is equal to $a$, let $a'$ denote the vertex of $\pi$ that immediately follows $a$ clockwise. Otherwise call $a'=a$. Observe that $\{a',b\}$ is an edge of $T$. It turns out that $U\contract{x}$ belongs to $\Omega(T\contract{x},a',b)$. Indeed, the deletion of $x$ cannot decrease the number of interior edges of $U$ incident to the other vertices of $\pi$ on the left of $(a,b)$. In addition, if this deletion makes an edge of $T$ and a previously distinct edge of $U$ identical, then the resulting edge must be incident to the vertex $y$ of $\pi$ that follows $x$ clockwise. However, as no interior edge of $U\contract{x}$ is incident to two vertices of $\pi$ on the left of $(a',b)$, and as $\{a',b\}$ is an edge of $T\contract{x}$, an interior edge common to $T\contract{x}$ and $U\contract{x}$ cannot be incident to a vertex of $\pi$ on the left of $\{a,b\}$. As a consequence, $y$ coincides with $b$ in this case, and all the interior edges shared by $T\contract{x}$ and $U\contract{x}$ must be incident to $a'$ or to $b$.

Now observe that there are still exactly $m$ vertices of $\pi$ on the right of $(a',b)$ that are adjacent to $a'$ or to $b$ by an edge of $T\contract{x}$ and to a vertex of $\pi$ that is not on the right of $(a',b)$ by an edge of $U\contract{x}$. Therefore, by induction,
\begin{equation}\label{Thm.far.1.eq.2}
d(T\contract{x},U\contract{x})\geq{(n-1)-m+(l-1)-5}\mbox{.}
\end{equation}

Combining (\ref{Thm.far.1.eq.1}) with (\ref{Thm.far.1.eq.2}) completes the proof.
\end{proof}

The bound provided by Theorem \ref{Thm.far.1} depends on the choice of two vertices. In order to choose these vertices in such a way that the resulting bound is as good as possible, the following notion is needed. Consider a triangulation $T$ of a polygon $\pi$ with $n$ vertices. A \emph{clockwise-oriented triangle of $T$} is an ordered triple $(a,b,c)$ of vertices of $\pi$ such that $c$ is on the right of $(a,b)$ and $\{a,b\}$, $\{b,c\}$, and $\{a,c\}$ are edges of $T$. Recall that, in this case, the lengths of $(a,b)$, $(b,c)$, and $(c,a)$ must sum to $n$. A clockwise-oriented triangle $(a,b,c)$ of $T$ is called \emph{central} when $(a,b)$, $(b,c)$, and $(c,a)$ all have length at most $n/2$. 

\begin{prop}\label{Prop.far.2}
Consider a triangulation $T$ of a convex polygon. At least one of the clockwise-oriented triangles of $T$ is central.
\end{prop}
\begin{proof}
Consider a triangulation $T$ of a polygon $\pi$ with $n$ vertices. If $n$ is equal to $3$, then $\pi$ admits a unique triangulation whose only triangle immediately provides a central clockwise-oriented triangle. Assume that $n$ is at least $4$. Let $(a,b,c)$ be a clockwise-oriented triangle of $T$ whose longest clockwise-oriented edge, say $(a,b)$, is the shortest possible among all the clockwise-oriented triangles of $T$.

It turns out that $(a,b,c)$ must be central. Indeed, the length of $(a,b)$ must be greater than $1$ because $n\geq4$. Hence, there exists a vertex $c'$ of $\pi$ on the right of $(b,a)$ such that $(b,a,c')$ is a clockwise-oriented triangle of $T$. Now assume for contradiction, that $(a,b)$ has length greater than $n/2$. Since the lengths of $(a,b)$ and $(b,a)$ sum to $n$, the length of $(b,a)$ is less than $n/2$. Further observe that the lengths of $(a,c')$ and $(c',b)$ sum to that of $(a,b)$. These two lengths are therefore also less than the length of $(a,b)$. This contradicts the assumption that the longest clockwise-oriented edge of $(a,b,c)$ is the shortest possible.
\end{proof}

Note that a triangulation of a convex polygon $\pi$ with $n$ vertices may have two distinct central clockwise-oriented triangles, up to reordering the vertices of these triangles. In this case, the two central triangles share an edge of length $n/2$, that splits $\pi$ into two polygons with the same number of vertices. In particular, $n$ must then be even and the triangulation cannot admit a third central triangle.

The following can now be stated and proven.

\begin{lem}\label{Lem.far.1}
Let $T$ be a triangulation of a polygon $\pi$. Consider a central clockwise-oriented triangle $(a,b,c)$ of $T$. Denote by $l$ the smallest length among those of $(a,b)$, $(b,c)$, and $(c,a)$. There exists a triangulation $U$ of $\pi$ such that
$$
d(T,U)\geq{n+l-6}\mbox{.}
$$
\end{lem}
\begin{proof}
Let $l_a$, $l_b$ and $l_c$ denote the lengths of $(a,b)$, $(b,c)$, and $(c,a)$, respectively. If $(x,y)$ is equal to $(a,b)$, to $(b,c)$, or to $(c,a)$, further call $m_x$ the number of vertices on the right of $(x,y)$ adjacent to $x$ or to $y$ by an interior edge of $T$. The sum $m_a+m_b+m_c-3$ counts the number of interior edges of $T$ incident to $a$, to $b$, or to $c$ (the subtracted $3$ stands for $\{a,b\}$, $\{b,c\}$, and $\{a,c\}$ being counted twice each). This sum cannot exceed the number of interior edges of $T$. In other words, 
$$
m_a+m_b+m_c\leq{n}\mbox{.}
$$

As $l_a$, $l_b$, and $l_c$ sum to $n$, it follows that
$$
2l_a+2l_b+2l_c\leq{3n-(m_a+m_b+m_c)}\mbox{.}
$$

Hence, there must be a vertex $x$ among $a$, $b$, and $c$ satisfying
\begin{equation}\label{Lem.far.1.eq.1}
2l_x\leq{n-m_x}
\end{equation}

It will be assumed, without loss of generality, that inequality (\ref{Lem.far.1.eq.1}) holds with $x=a$. Observe that the number of vertices of $\pi$ on the right of $(a,b)$ that are not adjacent to $a$ or to $b$ by an interior edge of $T$ is exactly $n-l_a-m_a-1$. Therefore, by (\ref{Lem.far.1.eq.1}), there are at least $l_a-1$ such vertices. Consider a polygon $\pi'$ whose vertices are $a$ and $b$ together with all the vertices of $\pi$ on the left of $(a,b)$, and exactly $l_a-1$ vertices on the right of $(a,b)$ that are not adjacent to $a$ or to $b$ by an interior edge of $T$.
\begin{figure}
\begin{centering}
\includegraphics{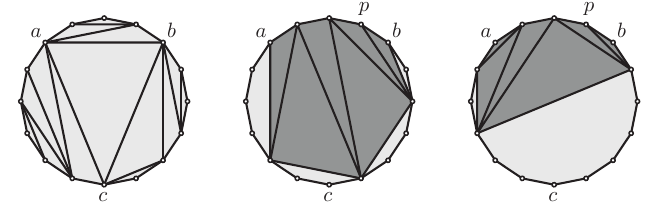}
\caption{An example of the construction in the proof of Lemma \ref{Lem.far.1}.
}\label{Far.fig.1}
\end{centering}
\end{figure}
Further consider a zigzag triangulation $Z$ of $\pi'$ with its ears in $a$ and in $b$: the interior edges of this triangulation form a simple path that alternates between left and right turns. There are two possible orientations for $Z$. Here, the chosen orientation is the one such that the vertex $p$ of $\pi$ preceding $b$ clockwise is incident to a single interior edge of $Z$. All the other vertices of $\pi$ on the left of $(a,b)$ are incident to exactly two interior edges of $Z$. This construction is illustrated in Fig. \ref{Far.fig.1} for two possible polygons $\pi'$ (colored dark grey), when $T$ and $(a,b,c)$ are as shown on the left of the figure. By construction, none of the (boundary or interior) edges of $Z$ are interior edges of $T$. Now complete $Z$ into a triangulation $U$ of $\pi$ by adding edges, in such a way that the number of interior edges shared by $T$ and $U$ is as small as possible. It turns out that $T$ and $U$ cannot share an interior edge. Indeed, otherwise, flipping such an edge within $U$ would result in a triangulation that contains all the edges of $Z$, and that shares one edge less with $T$. 

Recall that all the vertices of $\pi$ on the left of $(a,b)$ are incident to two interior edges of $Z$ except for $p$. In addition, note that $\{b,p\}$ cannot be incident to the same triangle in $T$ and in $U$ because the unique interior edge of $Z$ incident to $p$ crosses $\{a,b\}$. As a consequence, there must be at least one flip incident to $\{b,p\}$ along any path between $T$ and $U$. It then follows from Lemma \ref{Lem.intro.0} that:
\begin{equation}\label{Lem.far.1.eq.2}
d(T,U)\geq{d(T\contract{p},U\contract{p})+1}\mbox{.}
\end{equation}
 
Now observe that $(a,b)$, as an oriented edge on the polygon whose vertex set is obtained by removing $p$ from the vertex set of $\pi$, has length at most $\lceil{n/2}\rceil-1$. By construction, $U\contract{p}$ still does not share an interior edge with $T\contract{p}$, and all the vertices of $U\contract{p}$ on the left of $(a,b)$ are now incident to two interior edges of $U\contract{p}$. Hence, $U\contract{p}$ belongs to $\Omega(T\contract{p},a,b)$. Moreover, no vertex of $\pi$ on the right of $(a,b)$ is both adjacent to $a$ or to $b$ by an interior edge of $T\contract{p}$ and to a vertex that is not on the right of $(a,b)$ by an interior edge of $U\contract{p}$. In other words $m=0$ in the statement of Theorem~\ref{Thm.far.1} and this theorem yields
\begin{equation}\label{Lem.far.1.eq.3}
d(T\contract{p},U\contract{p})\geq{n-1+l_a-1-5}\mbox{.}
\end{equation}

Combining (\ref{Lem.far.1.eq.2}) with (\ref{Lem.far.1.eq.3}) and observing that $l_a\geq{l}$ completes the proof.
\end{proof}

When the central triangle in the statement of Lemma \ref{Lem.far.1} has a very short edge, the provided inequality is weak. In this case, the largest number of interior edges of a triangulation incident to any of its vertices comes into play. The following result takes this number into account.

\begin{lem}\label{Lem.far.2}
Consider a triangulation $T$ of a polygon $\pi$ with $n$ vertices. Let $(a,b,c)$ be a central clockwise-oriented triangle of $T$ and $l$ be the smallest length among those of $(a,b)$, $(b,c)$, and $(c,a)$. If every vertex of $\pi$ is incident to at most $n-3-k$ interior edges of $T$, then there exists a triangulation $U$ of $\pi$ such that
$$
d(T,U)\geq{n+\frac{k-9}{2}-l}\mbox{.}
$$
\end{lem}
\begin{proof}
Assume without loss of generality that $(c,a)$ has length $l$. Let $l_a$ and $l_b$ denote the lengths of $(a,b)$ and $(b,c)$, respectively. If $(x,y)$ is equal to $(a,b)$ or to $(b,c)$ further call $m_x$ the number of vertices on the right of $(x,y)$ adjacent to $x$ or to $y$ by an interior edge of $T$. Call $d$ the number of interior edges of $T$ incident to $b$. Note that $m_a+m_b-d-1$ counts the number of interior edges of $T$ incident to $a$ or to $c$ on one end, and to a vertex of $\pi$ that is not right of $(c,a)$ on the other. Again, $m_a+m_b-d$ needs to be subtracted by $1$ because $\{a,c\}$ would otherwise be counted twice. As there are at most $l-1$ edges whose two vertices are not on the right of $(c,a)$, the following inequality holds: 
\begin{equation}\label{Lem.far.2.eq.1.5}
m_a+m_b-d-1\leq{l-1}\mbox{.}
\end{equation}

Choose $k$ such that all the vertices of $\pi$ are incident to at most $n-3-k$ interior edges of $T$. Bounding $d$ accordingly and writing $n=l_a+l_b+l$, (\ref{Lem.far.2.eq.1.5}) yields
$$
l_a-m_a+l_b-m_b\geq{k-2l+3}\mbox{.}
$$

Therefore, the following holds with $x$ equal to $a$ or $b$:
\begin{equation}\label{Lem.far.2.eq.1}
l_x-m_x\geq\frac{k+3}{2}-l\mbox{.}
\end{equation}

It will be assumed without loss of generality that (\ref{Lem.far.2.eq.1}) holds with $x=a$. Consider a polygon $\pi'$ whose vertices are $a$ and $b$ together with all the vertices of $\pi$ on the left of $(a,b)$ and exactly $l_a-1$ vertices of $\pi$ on the right of $(a,b)$.
\begin{figure}[b]
\begin{centering}
\includegraphics{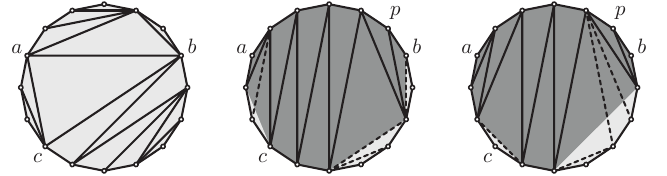}
\caption{An example of the construction in the proof of Lemma \ref{Lem.far.2}
}\label{Far.fig.2}
\end{centering}
\end{figure}
As in the proof of Lemma \ref{Lem.far.1}, consider a zigzag triangulation $Z$ of $\pi'$ with its ears in $a$ and in $b$. Observe that none of the interior edges of $Z$ belong to $T$ because they all cross $\{a,b\}$. Now complete the set of the interior edges of $Z$ into a triangulation $U$ of $\pi$ by adding edges, in such a way that the number of interior edges shared by $T$ and $U$ is as small as possible. Again, $T$ and $U$ cannot share an interior edge. Otherwise, flipping such an edge within $U$ would result in a triangulation that contains all the interior edges of $Z$, and that shares one less edge with $T$. Two triangulations $U$ built this way are shown in Fig. \ref{Far.fig.2} when $T$ is the triangulation represented on the left of the figure. Note that $\pi'$ is colored dark grey in each case, and that the interior edges of $U$ that are not interior edges of  $Z$ are dashed. As can be seen, the boundary edges of $\pi'$ do not necessarily belong to $U$.

As in the proof of Lemma \ref{Lem.far.1}, it can be assumed that the only vertex of $\pi$ on the left of $(a,b)$ that is possibly incident to less than two interior edges of $U$ is the vertex that immediately precedes $b$ clockwise. Call this vertex $p$ and observe that $\{b,p\}$ cannot be incident to the same triangle in $T$ and in $U$ because the interior edge of $Z$ incident to $p$ crosses $\{a,b\}$. As a consequence, there must be at least one flip incident to $\{b,p\}$ along any path between $T$ and $U$ and, by Lemma \ref{Lem.intro.0},
\begin{equation}\label{Lem.far.2.eq.2}
d(T,U)\geq{d(T\contract{p},U\contract{p})+1}\mbox{.}
\end{equation}
 
Again, the length of $(a,b)$, as an oriented edge on the polygon whose vertex set is obtained by removing $p$ from the vertex set of $\pi$, is at most $\lceil{n/2}\rceil-1$. Further note that the only interior edges common to $T\contract{p}$ and $U\contract{p}$ possibly created by the deletion of $p$ must be incident to $b$. Moreover, all the vertices of $U\contract{p}$ on the left of $(a,b)$ are incident to at least two (and this time possibly more than two) interior edges of $U\contract{p}$. Hence, $U\contract{p}$ belongs to $\Omega(T\contract{p},a,b)$. Moreover, the number of vertices of $\pi$ on the right of $(a,b)$ that are adjacent to $a$ or to $b$ by an interior edge of $T\contract{p}$ and to a vertex that is not on the right of $(a,b)$ by an interior edge of $U\contract{p}$ is at most $m_a$. Hence, Theorem \ref{Thm.far.1} yields:
$$
d(T\contract{p},U\contract{p})\geq{n-1-m_a+l_a-1-5}
$$

Here, $l_a-m_a$ can be bounded below using (\ref{Lem.far.2.eq.1}), implying
$$
d(T\contract{p},U\contract{p})\geq{n+\frac{k-11}{2}-l}
$$

Combining this inequality with (\ref{Lem.far.2.eq.2}) completes the proof.
\end{proof}

Lemmas \ref{Lem.far.1} and \ref{Lem.far.2} provide the following.

\begin{thrm}\label{Thm.far.2}
Consider a triangulation $T$ of a polygon $\pi$ with $n$ vertices. If all the vertices of $T$ are incident to at most $n-3-k$ interior edges of $T$, then $T$ has eccentricity at least $n+(k-21)/4$ in the flip-graph of $\pi$.
\end{thrm}
\begin{proof}
According to Proposition \ref{Prop.far.2}, $T$ admits a central clockwise-oriented triangle $(a,b,c)$. Let $l$ be the smallest length among those of $(a,b)$, $(b,c)$, and $(c,a)$. If $l$ is greater than, or equal to $(k+3)/4$, then Lemma \ref{Lem.far.1} provides the desired result, and if $l$ is less than $(k+3)/4$, then Lemma \ref{Lem.far.2} does. 
\end{proof}

Recall that Theorem \ref{MPThm.1} provides a lower bound on the eccentricity of a triangulation $T$ on the condition that some triangulation is compatible with a shelling of $T$ at a vertex $v$. By Lemma~\ref{Lem.sec1.2}, such a triangulation exists when $v$ is incident to more than half of the interior edges of $T$. In particular, Theorem \ref{MPThm.1} and Lemma \ref{Lem.sec1.2} do not provide any lower bound on the eccentricity of a triangulation when all of its vertices are incident to at most half of its interior edges. In contrast, the lower bound on eccentricities given by Theorem \ref{Thm.far.2}, while weaker, is valid for any triangulation. Moreover, it is large enough to allow for a precise characterization of the triangulations with a given, small enough eccentricity in the flip-graph of a convex polygon. 
Such a characterization is provided by Theorem \ref{Thm.Main.2}.

\begin{proof}[Proof of Theorem \ref{Thm.Main.2}]
Consider an integer $k$ such that $0\leq{k}\leq{n/8-5/2}$ and a triangulation $T$ of a polygon $\pi$ with $n$ vertices. First assume that some vertex of $\pi$ is incident to $n-3-k$ interior edges of $T$. Since $k\leq{n/2-2}$, it follows from Theorem \ref{Thm.Main.1} that the eccentricity of $T$ in the flip-graph of $\pi$ is $n-3+k$.

Now assume that the eccentricity of $T$ in the flip-graph of $\pi$ is $n-3+k$. Since $k\leq{n/8-5/2}$, this eccentricity is bounded above by
$$
\frac{9}{8}n-\frac{11}{2}\mbox{.}
$$

Let $n-3-l$ be the largest number of interior edges of $T$ incident to some of its vertices, and call $v$ such a vertex. By Theorem \ref{Thm.far.2}, the eccentricity of $T$ in the flip-graph of $\pi$ is at least $n+(l-21)/4$. As a consequence,
$$
n+\frac{l-21}{4}\leq\frac{9}{8}n-\frac{11}{2}\mbox{.}
$$

It follows that $l\leq(n-2)/2$. According to Theorem \ref{Thm.Main.1}, the eccentricity of $T$ in the flip-graph of $\pi$ is then $n-3+l$. This proves that $k$ is equal to $l$ and, therefore, that $v$ is incident to $n-3-k$ interior edges of $T$.
\end{proof}

\section{Concluding remarks and questions}\label{Sec.conc}

Let $\pi$ be a polygon with $n$ vertices. Further consider an integer $k$. According to Theorem \ref{Thm.Main.1}, all the triangulations of $\pi$ with $n-3-k$ interior edges incident to a vertex have the same eccentricity in the flip-graph of $\pi$ when $k\leq{n/2-2}$. The expression of this common eccentricity, $n-3+k$, is remarkably symmetric with the incidence number $n-3-k$. This seemingly unexpected result invites the question of whether something similar holds for the triangulations of $\pi$ whose largest number of interior edges incident to the same vertex is exactly $n-3-k$ with $k>n/2-2$. Note that the two cases are very different. When $k\leq{n/2-2}$, some vertex of $\pi$ is incident to more than half of the interior edges of $T$, and if the inequality is strict, then this vertex is necessarily unique. When $k>n/2-2$, there can be several (possibly many) vertices incident to exactly $n-3-k$ interior edges of of $T$.

It turns out that Theorem \ref{Thm.Main.1} does not extend to the latter case.

\begin{thrm}\label{Thm.rem.1}
Consider a polygon $\pi$ with $n$ vertices and an integer $k$ such that $n/2-2<k\leq{n-5}$. There exists a triangulation of $\pi$ whose largest number of interior edges incident to a vertex is exactly $n-3-k$, but whose eccentricity in the flip-graph of $\pi$ is at most $n-4+k$.
\end{thrm}
\begin{proof}
Call
$$
l=\left\lceil\frac{n-3}{n-3-k}\right\rceil+1\mbox{.}
$$

Since $k$ is greater than $n/2-2$ and at most $n-5$, one obtains
$$
3\leq{l}\leq\frac{n}{2}\mbox{.}
$$

In particular, there exists an oriented edge $(a,b)$ of length $l$ on $\pi$. Moreover, one can pick two vertices $x$ and $y$ of $\pi$ on the left of $(a,b)$ such that $\{x,y\}$ is an edge of $\pi$. Consider a triangulation $T$ of $\pi$ all of whose interior edges cross $\{a,b\}$. It follows from the choice of $l$ that the interior edges of $T$ can be placed in such a way that $x$ and $y$ are both incident to exactly $n-3-k$ interior edges of $T$, while all the other vertices of $\pi$ are incident to at most $n-3-k$ interior edges of $T$.

Consider a triangulation $U$ of $\pi$. Since $\{x,y\}$ is an edge of $\pi$, $U$ cannot have an ear in $x$ and an ear in $y$. Assume without loss of generality that $x$ is incident to at least one interior edge of $U$. In this case, $U$ can be transformed into the comb all of whose interior edges are incident to $x$ by a sequence of at most $n-4$ flips. The same triangulation can be reached from $T$ by exactly $k$ flips because $x$ is incident to exactly $n-3-k$ interior edges of $T$. As a consequence the distance of $T$ and $U$ in the flip-graph of $\pi$ is at most $n-4+k$, as desired.
\end{proof}

Theorem \ref{Thm.rem.1} further admits two interesting consequences. The first of these consequences is that Theorem \ref{Thm.Main.2} does not extend to all values of $k$.

\begin{cor}\label{Cor.rem.1}
Let $\pi$ be a polygon with $n$ vertices. There exists an integer $k$ such that $n/8-5/2<k\leq{n/2-2}$ and a triangulation $T$ of $\pi$ whose vertices are all incident to at most $n-4-k$ interior edges, but whose eccentricity in the flip-graph of $\pi$ is exactly $n-3+k$.
\end{cor}
\begin{proof}
Invoking Theorem \ref{Thm.rem.1} with $k=\lfloor{n/2}\rfloor-1$, one obtains a triangulation $T$ of $\pi$ whose largest number of interior edges incident to a vertex is exactly $\lceil{n/2}\rceil-2$, and whose eccentricity in the flip-graph of $\pi$ is at most $\lfloor{3n/2}\rfloor-5$. Hence
$$
k\leq\left\lfloor\frac{n}{2}\right\rfloor-2\mbox{,}
$$
where $k$ denotes the number obtained by subtracting $n-3$ from the eccentricity of $T$ in the flip-graph of $\pi$. This inequality can be rewritten as
$$
\left\lceil\frac{n}{2}\right\rceil-2\leq{n-4-k}\mbox{.}
$$

In other words, the vertices of $T$ are all incident to at most $n-4-k$ edges. Finally, by Theorem \ref{Thm.Main.2}, $k$ must be greater than $n/8-5/2$.
\end{proof}

The other consequence is that the largest number of interior edges incident to a vertex can be the same for two triangulations of a polygon while their eccentricities in the flip-graph of that polygon are distinct.

\begin{cor}
Let $\pi$ be a polygon with $n$ vertices. If $n$ is greater than $12$, then there exist two triangulations of $\pi$ whose largest number of interior edges incident to one of their vertices is four, while their eccentricities in the flip-graph of $\pi$ are respectively exactly $2n-10$ and at most $2n-11$.
\end{cor}
\begin{proof}
By Theorem \ref{Thm.rem.1}, there exists a triangulation of $\pi$ whose largest number of interior edges incident to a vertex is exactly four, but whose eccentricity in the flip-graph of $\pi$ is at most $2n-11$. It is proven in \cite{Pournin2014}, though that, when $n>12$, there exist triangulations of $\pi$ whose largest number of interior edges incident to a vertex is four, but whose eccentricity in the flip-graph of $\pi$ is $2n-10$.
\end{proof}

In view of the above results, it is natural to ask the following questions.

\begin{question}
Let $\pi$ be a polygon with $n$ vertices. What are the values of $k$ such that a triangulation of $\pi$ has eccentricity $n-3+k$ in the flip-graph of $\pi$ if and only if the largest number of its interior edges incident to a vertex is $n-3-k$?
\end{question}

\begin{question}\label{Ques.rem.2}
Let $\pi$ be a polygon with $n$ vertices and $k$ an integer such that $n/2-2<k\leq{n-5}$. Consider the triangulations of $\pi$ whose largest number of interior edges incident to a vertex is exactly $n-3-k$. What are all the possible values of the eccentricity in the flip-graph of $\pi$ of these triangulations ?
\end{question}

Proving that there is a triangulation of $\pi$ with eccentricity $n-3+k$ whenever $n/2-2<k<n-7$ would provide an interesting partial answer to Question \ref{Ques.rem.2}. Note that, when $k$ is equal to $n-7$, such triangulations are given in \cite{Pournin2014}. Further note that, when $k$ is equal to $n-6$ or to $n-5$, the upper bound of $n-3+k$ on the eccentricity of a triangulation of $\pi$ provided by Proposition \ref{Intro.prop.1} cannot be sharp because it is greater than the diameter of the flip-graph of $\pi$. In these cases, Question \ref{Ques.rem.2} is interesting nonetheless. In particular, not much is known about the eccentricity of zigzag triangulations in the flip-graph of a convex polygon.

\begin{acknowledgement}I thank Thibault Manneville for many discussions on the subject of flip distances, that inspired the results presented in this article. I also thank the two anonymous referees of this article for providing a number of very useful comments, that helped improving the exposition of the results.
\end{acknowledgement}

\end{document}